\documentclass[12pt,a4paper]{amsart}
\usepackage{graphicx}
\usepackage{float,epsfig}
\usepackage[russian,english]{babel}
\usepackage[cp866]{inputenc}
\usepackage[psamsfonts]{amssymb}
\usepackage{amsfonts}
\usepackage{amsmath}
\newtheorem{theorem}{Theorem}
\newtheorem{lemma}{Lemma}
\newtheorem{corollary}{Corollary}
\def\Re{\mathop{\rm Re}\nolimits}
\def\Im{\mathop{\rm Im}\nolimits}

\newcounter{minutes}
 \divide\time by 60
\newcounter{hours}
 \multiply\time by 60
\addtocounter{minutes}{-\time}

\begin{document}

\title[]{CONFORMAL MODULE OF THE EXTERIOR \\ OF TWO RECTILINEAR
SLITS}\thanks{The work of the first author was supported by the
Russian Foundation for Basic Research and the Government of the
Republic of Tatarstan, grant No~18-41-160003; the second author
was supported by the Russian Foundation for Basic Research, grant
No~17-01-00282. The third author expresses his thanks to the Kazan
Regional Scientific and Educational Mathematical Center for a
support during his stay at Kazan Federal University in
October-November 2018.}

\author[D.~Dautova]{D.~Dautova}
\address{Kazan Federal University,
         Kazan, Russia}
         \email[dautovadn@gmail.com]{dautovadn@gmail.com}
\author[S.~Nasyrov]{S.~Nasyrov}
\address{Kazan Federal University,
         Kazan, Russia}
 \email[]{semen.nasyrov@yandex.ru}

\author[M.~Vuorinen]{M.~Vuorinen}
\address{Department of Mathematics and Statistics, University of Turku,
         Turku, Finland}
\email[vuorinen@utu.fi]{vuorinen@utu.fi}






\maketitle

\begin{abstract}
We study moduli of planar ring domains whose complements are
linear segments and establish formulas for their moduli in terms
of the Weierstrass elliptic functions. Numerical tests are carried
out to illuminate our results.


\end{abstract}



\section{Introduction}\label{intr}

The Weierstrass and Jacobian elliptic and theta functions and the
Schwarz-Christoffel formula form the foundation for numerous
explicit formulas for conformal mappings (N.I. Akhiezer
\cite{akhiezer}, W.~Koppenfels, F.~Stallmann \cite{kop_sht}).
During the past thirty years many authors have studied numerical
implementation of conformal mappings. We refer the reader to the
bibliography of the monograph  N. Papamichael and  N.
Stylianopoulos \cite{ps}. In particular, the Schwarz-Christoffel
toolbox of T. Driscoll and  N. Trefethen \cite{dt} has become a
standard tool in the field. In a series of papers of T. DeLillo,
J. Pfaltzgraff, D. Crowdy and their coauthors have extended the
Schwarz-Christoffel method to certain cases of multiply connected
domains with polygonal boundary components
\cite{cr,delil2,delil1}.

In addition to the conformal mapping problem,  also the
computation of numerical values of conformal invariants is an
important issue in geometric function theory. Here one can often
use a conformal map onto a canonical domain so as to simplify the
computation. Therefore computation of conformal invariants has a
natural link to numerical conformal mapping. However, the so
called crowding phenomenon can create serious obstacles for
computation of conformal maps, e.g. when long rectangles are mapped onto the upper half space
\cite{ps}.

A basic conformal invariant is the module of a ring domain. A ring
domain $G$ can be conformally mapped onto an annulus $\{z \in
\mathbb{C}: q<|z|<1\}$ and its conformal  module and capacity are
defined as
$$
{\rm mod}\, G =( \log (q^{-1}))/(2\pi)\,, \quad {\rm cap}\, G = 2
\pi/\log (q^{-1}) \,.
$$
Therefore, ${\rm mod}\,G =1/{\rm cap}\, G$ and the computation of
${\rm mod}\, G$ can be reduced to the solution of the Dirichlet
problem for the Laplace equation and to the computation of the
$L^2$-norm of its gradient. This method was applied in
\cite{bsv,hrv1,hrv3} for the case of bounded ring domains.

Here we shall consider unbounded ring domains whose complementary
components are segments. We describe one-parametric families of
functions $f(z,t)$ each of which maps conformally  an annulus
$\{q<|\zeta|<1\}$ onto the exterior $G=G(t)$ of two disjoint
segments $A_1A_2$ and $A_3A_4$. Here $A_j=A_j(t)$, $1\le j \le 4$,
are some smooth functions and $q=q(t)$; $t$ is a real parameter.
Further we will denote  such domains by $G(A_1,A_2,A_3,A_4)$. It
is also assumed that the straight lines, containing the segments
$A_1A_2$ and $A_3A_4$, are fixed.

We note that one-parametric families of conformal mappings were
considered earlier. There is the well-known Loewner-Komatu
differential equation which is a generalization of the Loewner
equation to the doubly-connected case. The approach of Komatu was
developed by Goluzin~\cite{gol1} and others (see, e.g.
\cite{alex,gum1,gum2}).

We deduce a differential equation for $f(z,t)$ in the considered
case
 (Theorem~\ref{loew-ell}). In contrast to the
Loewner-Komatu equation, we do not assume that the family of the
images is monotonic as a function of the parameter $t$. As a
corollary, we obtain a system of ODEs to determine the behavior of
the accessory parameters, which are the preimages of the points
$A_j$, and the conformal module $m(t):={\rm mod}\, G(t) =(\log
(q(t))^{-1})/(2\pi)$. On the base of the system, we suggest an
approximate method for finding the accessory parameters and the
conformal module. We note that in our approach we use essentially
the Weierstrass elliptic functions.

Further we apply the obtained  results to investigate the behavior
of the conformal module in the case when one of the segments and
the length of the other one are fixed.

Now we briefly describe the structure of the paper. In
Section~\ref{prel} we give  some information on the Weierstrass
elliptic functions, moduli, and reduced moduli. In
Section~\ref{repr} we describe an integral representation of an
annulus onto the exterior of two slits $A_1A_2$ and $A_3A_4$
(Theorem~\ref{map_fix}). In contrast to the known representations
\cite{komatu}; (see also \cite{henrici3},  \cite{delil1},
\cite{delil2}), our representation is based on the Weierstrass
$\sigma$-functions. The representation contains some unknown
constants; they are called accessory parameters. In
Section~\ref{fam} we consider one-parametric families of such
functions $f(z,t)$ and deduce a differential equation for them
(Theorem~\ref{loew-ell}). As a corollary, we obtain a system of
ODEs for accessory parameters (Theorem~\ref{system}). We note that
to deduce the equations we use the approach developed earlier for
one-parametric families of rational functions \cite{nas_dokl} and
conformal mappings of complex tori \cite{nas_vuz,nas_petr}. In
Section~\ref{numeric} we give results of some numerical
calculation. In Section~\ref{extr} we study monotonicity of the
conformal module of the exterior of two slits when one segment is
fixed and the other one slides along a straight line and has a
fixed length.

Finally we should note that recently the capacity computation of
doubly connected domains with complicated boundary structure has
been studied for instance in \cite{hrv3}.

\section{Some preliminary results}\label{prel}

\textit{Elliptic functions.} First we recall some information
about elliptic functions (see, e.g., \cite{akhiezer,nist} and also
\cite{dlmf1,dlmf2}).

A meromorphic in the complex plane function is called elliptic if
it has periods $\omega_1$ and
$\omega_2$\footnotemark{*}\footnotetext{${}^*$ In contrast to
\cite{akhiezer}, we denote by $\omega_1$ and $\omega_2$ periods of
elliptic functions, not half-periods. The same remark concerns the
values $\eta_k$ defined by (\ref{periodzeta}).}, linearly
independent over $\mathbb{R}$. In the fundamental parallelogram
constructed by the vectors $\omega_1$ and $\omega_2$, every
nonconstant elliptic function takes each value the same number of
times; the number $r$ is called the order of the elliptic
function.

If $a_1,\ldots,a_r$ are zeroes of an elliptic functions of order
$r$ and $b_1,\ldots,b_r$ are its poles in the fundamental
parallelogram, then
$$
a_1+\ldots+a_r\equiv b_1+\ldots+b_r\quad (\textrm{mod}\,\Omega)
$$
where $\Omega$ is the lattice generated by  $\omega_1$ and
$\omega_2$. Further we will denote by $\omega$ an arbitrary
element of the lattice. We note that, by given lattice, the
generators $\omega_1$ and $\omega_2$ are not determined by a
unique way; we will further assume that
$\Im(\omega_2/\omega_1)>0$.

One of the main elliptic functions is the Weierstrass
$\mathfrak{P}$-function
\begin{equation*}\label{pw}
\mathfrak{P}(z)=\frac{1}{z^2}+\sum\nolimits'\left[\frac{1}{(z-\omega)^2}-\frac{1}{\omega^2}\right];
\end{equation*}
here the summation $\sum'$ is over all nonzero elements of the
lattice. The Weierstrass $\zeta$-function
\begin{equation}\label{zw}
\zeta(z)=\frac{1}{z}+\sum\nolimits'\left[\frac{1}{z-\omega}+\frac{1}{\omega}+\frac{z}{\omega^2}\right]
\end{equation}
has the properties: $\zeta'(z)=-\mathfrak{P}(z)$ and
\begin{equation}\label{periodzeta}
\zeta(z+\omega_k)=\zeta(z)+\eta_k,\quad k=1,2,
\end{equation}
where $\eta_k=2\zeta(\omega_k/2)$. In the fundamental
parallelogram it has a unique pole with residue~$1$. The numbers
$\eta_k$ and $\omega_k$ satisfy the equality
\begin{equation}\label{rel}
\omega_2\eta_1-\omega_1\eta_2=2\pi i.
\end{equation}
At last, we need the Weierstrass $\sigma$-function
\begin{equation}\label{sw}
\sigma(z)=z\prod\nolimits'\left\{\left(1-\frac{z}{\omega}\right)\exp\left(\frac{z}{\omega}+\frac{z^2}{2\omega^2}\right)\right\}.
\end{equation}
It is an odd entire function with the properties:
$$
\frac{\sigma'(z)}{\sigma(z)}=\zeta(z),\quad
\sigma(z+\omega)=\varepsilon \sigma(z)e^{\eta(z+\omega/2)}.
$$
Here $\eta=m\eta_1+n\eta_2$, if $\omega=m\omega_1+n\omega_2$.
Moreover, $\varepsilon=1$, if $\omega/2$ belongs to the lattice
$\Omega$, otherwise, $\varepsilon=-1$.

We recall the Weierstrass invariants $g_2$ and $g_3$:
\begin{equation*}
g_2=60\sum\nolimits'\frac{1}{(m\omega_1+n\omega_2)^4},\quad
g_3=140\sum\nolimits'\frac{1}{(m\omega_1+n\omega_2)^6}.
\end{equation*}

Elliptic functions depend  not only on the variable $z$ but also
on the lattice. Further we need explicit expressions for the
partial derivatives $\zeta(z)= \zeta(z;\omega_1,\omega_2)$ by the
periods $\omega_1$ and $\omega_2$ of the lattice. In
\cite{nas_vuz}   the following theorem is proved.

\begin{theorem}\label{defzeta}
The partial derivatives  of $\zeta(z)=\zeta(z;\omega_1,\omega_2)$
with respect to the periods $\omega_1$ and $\omega_2$ are equal to
\begin{equation*}\label{zetaom1}
\frac{\partial\zeta(z)}{\partial \omega_1}=\frac{1}{2\pi i}\left[
\frac{1}{2}\,\omega_2\mathfrak{P}'(z)+(\omega_2\zeta(z)-\eta_2z)\mathfrak{P}(z)+
\eta_2 \zeta(z)-(\omega_2 g_2/12)z\right],
\end{equation*}
\begin{equation*} \label{zetaom2}\frac{\partial\zeta(z)}{\partial
\omega_2}=-\frac{1}{2\pi i}\left[
\frac{1}{2}\,\omega_1\mathfrak{P}'(z)+(\omega_1\zeta(z)-\eta_1z)\mathfrak{P}(z)+
\eta_1 \zeta(z)-(\omega_1 g_2/12)z\right].
\end{equation*}
\end{theorem}

We will also need the Jacobi theta-function $\vartheta_1(z)$. For
given lattice $\Omega$, generated by $\omega_1$ and $\omega_2$,
let $\tau=\omega_2/\omega_1$, $\Im \tau>0$, and $q=e^{\pi i
\tau}$. Then, by definition (see, e.g. \cite[ch.~1,
sect.~3]{akhiezer}, \cite{dlmf1}),
\begin{equation}\label{theta1}
\vartheta_1(z)=\vartheta_1(z|\tau)=1+\sum_{n=0}^\infty
(-1)^nq^{(n+1/2)^2}\sin((2n+1)z).
\end{equation}
There is the following connection between $\sigma$-function and
$\vartheta_1(z)$ (see, e.g. \cite[ch.~4, sect. 19, formula
(1)]{akhiezer}):
\begin{equation}\label{connect}
\sigma(z)=\omega_1 \frac{e^{\frac{\eta_1
z^2}{2\omega_1}}\vartheta_1(z/\omega_1 )}{\vartheta'_1(0)}\,.
\end{equation}

\vskip 0.4 cm

\textit{Conformal moduli, reduced moduli, and capacities of
condensers.} Let $G$ be a ring domain  in the plane, i.e. a
doubly-connected domain with non-degenerate boundary components.
There is a conformal mapping $\psi:G\to A$ of $G$ onto an annulus
$A=\{q<|z|<1\}$ (see, e.g., \cite{gol}). The value $q$ does not
depend on the choice of $\psi$. We call $$
\text{mod}\,G=\frac{1}{2\pi}\,\log (q^{-1})
$$
the conformal module of $G$. It is conformal invariant and plays
an important role in the theory of conformal and quasiconformal
mappings.

Let $G$ be a ring domain in the plane with complementary
components $C_1$ and $C_2\,,$ and let $K$ be the condenser with
plates $C_1$ and $C_2\,$ and with field $G\,.$ We recall that $$
\text{cap}\,K=\inf\limits_u \int\!\!\!\int|\nabla u|^2 dxdy
$$
where the infimum is taken over all smooth functions $u$ such that
$u=0$ on $C_1$ and $u=1$ on $C_2$. We will define
$\text{cap}\,G:=\text{cap}\,K$ and call $\text{cap}\,G$ the
conformal capacity of the ring domain~$G\,.$

Let $D$ be a simply connected domain with non-degenerate boundary
and $z_0\in D$. For sufficiently small $\varepsilon$ consider the
condenser defined by $D\setminus B_{\varepsilon}(z_0)$; here
$B_{\varepsilon}(z_0)$ is the disk of radius $\varepsilon$
centered at the point $z_0$. Denote by $K_\varepsilon$ its
capacity.  Then there exists the limit
$$ r(D,z_0):=\lim_{\varepsilon\to 0+}(K_\varepsilon+(1/(2\pi))\log
\varepsilon)
$$
which is called the reduced module of $D$ at the point $z_0$
\cite[Section 2.4]{dubook}, \cite{garnett}.

\section{Integral representation}\label{repr}

\begin{figure}[ht] \centering
\includegraphics[width=4.5 in,%
]{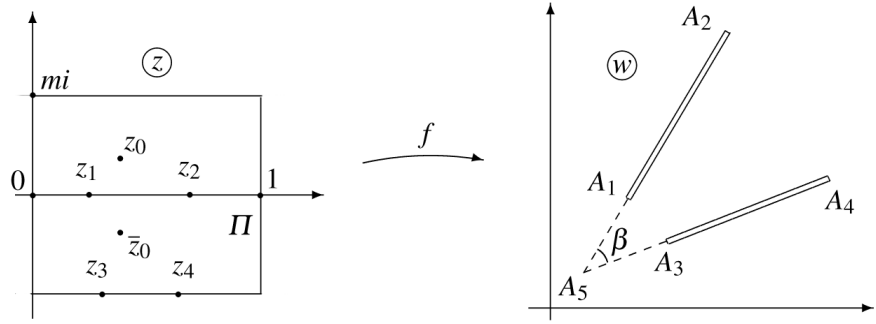}\ \caption{Conformal mapping of the rectangle $\Pi$
with identified vertical sides  onto $G(A_1,A_2,A_3,A_4)$.}
\label{confo}
\end{figure}

Consider a conformal mapping $g$ of an annulus $\{q<|\zeta|<1\}$
onto the exterior $G=G(A_1,A_2,A_3,A_4)$ of two disjoint
rectilinear slits $A_1A_2$ and $A_3A_4$ in the $w$-plane. With the
help of the exponential map $z\mapsto \zeta=\exp(2\pi i z)$ we can
consider the map $f:=g(2\pi i z)$ from the horizontal strip
$$S:=\{-m<\Im z<0\},\quad m=\frac{1}{2\pi}\,\log (q^{-1}),$$ onto $G$. It
maps conformally the rectangle $\Pi=\{0<\Re z<1, \,-m<\Im z<0\}$
with identified vertical sides onto $G$ (Fig.~\ref{confo}).  The
value $m$ is the conformal module of $G$. It is evident that $f$
has a unique pole in $\Pi$.

We will find an integral representation for the  conformal mapping
$f$ of Schwarz--Christoffel type using the Weierstrass
$\sigma$-function. We should note that analogs of the
Schwarz--Christoffel integral for doubly-connected domains were
obtained earlier in \cite{komatu}; it is based on
$\theta$-functions (see also \cite{delil1,delil2,henrici3}).

Using the Riemann--Schwarz reflection principle, we can extend $f$
to $\mathbb{C}$ as a meromorphic function. We see that the
function $h(z)=f''(z)/f'(z)$ is doubly periodic in $\mathbb{C}$
with periods $\omega_1=1$ and $\omega_2=2mi$. Consider $h$ in the
double rectangle $\widetilde{\Pi}=\{0<\Re z<1, \,-m<\Im z<m\}$; it
is its fundamental parallelogram.  Here the function $h$ has only
polar singularities at points $z_k$, $1\le k\le 4$, corresponding
to the endpoints $A_k$ of the slits, and also at two distinct
points, $z_0$ and $\overline{z}_0$, where $f$ has poles. For
definiteness, we assume that $y_0:=\Im z_0>0$. The residues of $h$
are known, therefore, we can express it with the help of the
Weierstrass zeta-function:
\begin{equation}\label{h}
h(z)=\gamma+\sum_{k=1}^4\zeta(z-z_k)-2\zeta(z-z_0)-2\zeta(z-\overline{z}_0)
\end{equation}
where $\gamma$ is a constant. (Here and further, unless otherwise
specified,  we assume that $\zeta(z)$ and other elliptic functions
have periods $\omega_1=1$ and $\omega_2=2mi$.)

From (\ref{h}) we have
\begin{equation*}
\log f'(z)= \gamma z+\log
C+\sum_{k=1}^4\log\sigma(z-z_k)-2\log\sigma(z-z_0)-2\log\sigma(z-\overline{z}_0),
\end{equation*}
\begin{equation*}\label{f_prime}
f'(z)= C e^{\gamma
z}\frac{\prod_{k=1}^4\sigma(z-z_k)}{\sigma^2(z-z_0)\sigma^2(z-\overline{z}_0)},
\end{equation*}
\begin{equation}\label{f}
f(z)= C\int_{0}^z e^{\gamma
\xi}\frac{\prod_{k=1}^4\sigma(\xi-z_k)}{\sigma^2(\xi-z_0)\sigma^2(\xi-\overline{z}_0)}\,d\xi+C_1
\end{equation}
where $\sigma(z)$ is the Weierstrass sigma-function, $C\neq0$ and
$C_1$ are complex constants.

The residue of $f'(z)$ at $z_0$ must vanish, therefore,
\begin{equation*}\label{res}
\gamma
+\sum_{k=1}^4\zeta(z_0-z_k)-2\zeta(z_0-\overline{z}_0)=0\,.
\end{equation*}

The $\sigma$-function satisfies (\ref{periodzeta}).
 Because $f'(z)$ must be periodic
with period $\omega_1=1$, we have
\begin{multline*}
f'(z+1)= C e^{\gamma(
z+1)}\frac{\prod_{k=1}^4\sigma(z-z_k+1)}{\sigma^2(z-z_0+1)\sigma^2(z-\overline{z}_0+1)}\,\\=\frac{Ce^{\gamma}
e^{\gamma
z}\prod_{k=1}^4e^{\eta_1(z-z_k+1/2)}\sigma(z-z_k)}{e^{2\eta_1(z-z_0+1/2)}\sigma^2(z-z_0)e^{2\eta_1(z-\overline{z}_0+1/2)}\sigma^2(z-\overline{z}_0)}\,
=e^{\gamma+\eta_1(2z_0+2\overline{z}_0-\sum_{k=1}^4z_k)}f'(z).
\end{multline*}
Consequently,
\begin{equation}\label{1}
\gamma+\eta_1\left(2z_0+2\overline{z}_0-\sum_{k=1}^4z_k\right)\equiv
0\quad (\mbox{\rm mod}\, 2\pi i).
\end{equation}

In a similar way, we have
\begin{equation*}
f'(z+\omega_2)=e^{\gamma\omega_2+\eta_2(2z_0+2\overline{z}_0-\sum_{k=1}^4z_k)}f'(z).
\end{equation*}
Since $\arg f'(z+\omega_2)-\arg f'(z)=2\beta$ where $\beta$ is the
angle between the segments $A_1A_2$ and $A_3A_4$, we have
\begin{equation}\label{2}
\gamma\omega_2+\eta_2\left(2z_0+2\overline{z}_0-\sum_{k=1}^4z_k\right)=2\beta
i \quad (\mbox{\rm mod}\, 2\pi i).
\end{equation}
Now we will specify the position of the points $z_k$. We will
assume that $z_1$ and $z_2$ lie on the real axis, and $z_3$ and
$z_4$ are on the lower side of  $\Pi$ (Fig.~\ref{confo}). Because
$z_k$ can be chosen up to the values $k \omega_1+n\omega_2$, $k$,
$n\in \mathbb{N}$, for convenience, we shift $z_3$ by
$\omega_2=2mi$ and assume that
\begin{equation}\label{xk}
z_1=x_1,\ z_2=x_2,\  z_3=x_3+i m,\ z_4=x_4-i m,
\end{equation}
where $x_k$ are real numbers.

Denote
\begin{equation}\label{a1}
a=\sum_{k=1}^4z_k-2z_0-2\overline{z}_0.
\end{equation}

Taking into account (\ref{xk}), we see that $a$ is real. We write
(\ref{1}) and (\ref{2}) in the form
\begin{equation}\label{sys}
\gamma- \eta_1 a =2\pi k i,\quad \omega_2\gamma-\eta_2 a=2\beta i
+2\pi n i,\quad k,n\in \mathbb{N}.
\end{equation}
Solving (\ref{sys}) as a system of linear equation with respect to
$\gamma$ and $a$ and taking into account that its determinant
equals $\omega_1\eta_2-\omega_2\eta_1=2\pi i$, we obtain
\begin{equation}\label{sys1}
\gamma=-k\eta_2 +(n+\beta/\pi)\eta_1, \quad
a=-k\omega_2+(n+\beta/\pi).
\end{equation}
Since $\omega_2$ is a purely imaginary number, from the second
equality in (\ref{sys1}) we deduce that $k=0$. We can change $x_3$
by entire values, therefore, we can assume that  $n=0$. So
(\ref{sys1}) has the form
\begin{equation}\label{sys2}
\gamma=\beta\eta_1/\pi, \quad a=\beta/\pi.
\end{equation}

Thus, from (\ref{a1}) and (\ref{sys2}) we have
\begin{equation*}
\sum_{k=1}^4x_k=4x_0+\beta/\pi,\quad x_0=\Re z_0.
\end{equation*}

Since every horizontal shift does not change the strip $S$, we can
assume that $x_0=0$.

Therefore, we establish the following theorem.

\begin{theorem}\label{map_fix}
The function, mapping the annulus $\{q<|\zeta|<1\}$ onto 
$G(A_1,A_2,A_3,A_4)$,  is $f(z)$ where $z=(2\pi i)^{-1}\log\zeta$
and $f$ is defined by (\ref{f}). In (\ref{f})
$\gamma=\beta\eta_1/\pi$, the points $z_k=x_k+i y_k$ correspond to
the endpoints $A_k$ of the slits and satisfy (\ref{xk}) with real
$x_k$ and $m=(1/(2\pi))\log (q^{-1})$, the point $z_0=i y_0$
matches to the infinity, $C\neq 0$ and $C_1$ are some complex
constants. Moreover, $\sum_{k=1}^4x_k=\beta/\pi$.
\end{theorem}

\section{One-parametric families}\label{fam}
The parametric method  for doubly connected domains was developed
by Komatu~\cite{komatu1} and Goluzin~\cite{gol1} (in details, see
\cite{alex}, ch.~5). In recent papers \cite{gum1,gum2} some new
results were obtained. Here we obtain an equation of Loewner type
using ideas of the papers \cite{nas_vuz,nas_dokl}.

Taking into account the integral representation (\ref{f}),
obtained in Theorem~\ref{map_fix}, we consider a smooth
one-parametric family of conformal mappings
\begin{equation}\label{family}
f(z,t)= c(t)\int_{0}^z e^{\gamma(t)
\xi}\frac{\prod_{k=1}^4\sigma(\xi-z_k(t))}{\sigma^2(\xi-z_0(t))\sigma^2(\xi-\overline{z_0(t)})}\,d\xi+c_1(t)
\end{equation}
Here $\sigma(z)=\sigma(z;1,\omega_2)$ where $\omega_2=2mi$,
$m=m(t)>0$. For a fixed $t$, $f(z,t)$ is periodic with period
$\omega_1\equiv 1$ and maps the half of the fundamental
parallelogram (rectangle) $\{0<\Re x<1, -m<\Im z<0\}$ onto the
exterior of two rectilinear slits. Without loss of generality we
may assume that one slit lies on the positive part of the real
axis and the second one is on the ray $\{\arg w=\beta\}$. (The
general case can be obtained by multiplying $c(t)$ by
$e^{i\theta}$; this means the rotation by the angle $\theta$.
Further, in some situations, we will use this remark.)

We note once more that the angle $0<\beta<2\pi$ does not depend on
$t$. Moreover,
\begin{equation*}
\Im z_1(t)=\Im z_2(t)=0,\quad \Im z_3(t)=-\Im z_4(t)=m(t), \quad
\Re z_0(t)=0,
\end{equation*}
therefore,
\begin{equation*}
z_1(t)=x_1(t),\  z_2(t)=x_2(t),\ z_3(t)=x_3(t)+im(t),\
z_4(t)=x_4(t)-im(t),\ z_0(t)=iy_0(t),
\end{equation*}
\begin{equation*}
x_1(t)<x_2(t)<x_1(t)+1,\ x_3(t)<x_4(t)<x_3(t)+1,\ 0\le y_0(t)\le
m,
\end{equation*}
\begin{equation*}
\gamma(t)=(\beta/\pi)\eta_1(t), \quad
\sum_{k=1}^4x_k(t)=\beta/\pi,
\end{equation*}
\begin{equation*}
\gamma(t)
+\sum_{k=1}^4\zeta(z_0(t)-z_k(t))-2\zeta(z_0(t)-\overline{z_0(t)})=0,
\end{equation*}

By the Riemann-Schwarz symmetry principle, we can extend $f(z,t)$
meromorphically  to the whole complex plane. It is evident that
the extension satisfies
\begin{equation}\label{quasiper}
f(z+1,t)=f(z,t),\quad f(z+\omega_2(t),t)=e^{2i\beta}f(z,t),
\end{equation}
Differentiating \eqref{quasiper} with respect to $t$ and $z$, we obtain
$$
\dot{f}(z+1,t)=\dot{f}(z,t), \quad
\dot{\omega}_2(t)f'(z+\omega_2(t),t)+\dot{f}(z+\omega_2(t),t)=e^{2i\beta}\dot{f}(z,t),
$$
\begin{equation*}
f'(z+1,t)=f'(z,t),\quad f'(z+\omega_2(t),t)=e^{2i\beta}f'(z,t).
\end{equation*}
Here and further the dot means differentiation  with respect to
the parameter $t$ and the prime is differentiation  with respect
to $z$. Thus, we have
$$
\frac{\dot{f}(z+\omega_k(t),t)}{f'(z+\omega_k(t),t)}+\dot{\omega}_k(t)=\frac{\dot{f}(z,t)}{f'(z,t)}.
$$

Consequently, the function $h(z,t):={\dot{f}(z,t)}/{f'(z,t)}$
satisfies
\begin{equation}\label{periodh}
h(z+\omega_k(t),t)-h(z,t)=-\dot{\omega}_k(t),\quad k=1,2,
\end{equation}
where $\dot{\omega}_1(t)\equiv 0$.

Now we write Taylor's expansion of $f(z,t)$ in a neighborhood of
$z_k(t)$:
\begin{equation}\label{taylor}
 f(z,t)=A_k(t)+\frac{D_k(t)}{2}\,(z-z_k(t))^2+\ldots,\quad
\end{equation}
{where} $D_k(t)=f''(z_k(t),t). $ We have
\begin{multline*}
f''(z,t)=c(t)e^{\gamma(t) z
}\,\frac{\prod\limits_{j=1}^4\sigma(z-z_j(t))}{\sigma^{2}(z-z_0(t))\sigma^{2}(z-\overline{z_0(t)})}\,\\
\times
\Bigl[\,\gamma(t)+\sum_{j=1}^4\zeta(z-z_j(t))-2\zeta(z-z_0(t))-2\zeta(z-\overline{z_0(t)})\Bigr],
\end{multline*}
therefore, as $z\to z_k(t)$, we obtain
\begin{equation}\label{dk}
D_k(t)= c(t)e^{\gamma(t) z_k(t)}\,\frac{\prod\limits_{j=1,\,j\neq
k}^4\sigma(z_k(t)-z_j(t))}{\sigma^{2}(z_k(t)-z_0(t))\sigma^{2}(z_k(t)-\overline{z_0(t))}}\,.
\end{equation}
From (\ref{taylor}) it follows that
\begin{equation}\label{fprime}
 f'(z,t)=D_k(t)(z-z_k(t))+\ldots,
\end{equation}
\begin{equation*}\label{(3.2)}
 \dot{f}(z,t)=\dot{A}_k(t)-\dot{z}_k(t)D_k(t)(z-z_k(t))+\ldots,
\end{equation*}
and, therefore,
\begin{equation*} \label{(3.3)}
 h(z,t)=\frac{\dot{f}(z,t)}{f'(z,t)}=\frac{\gamma_k(t)}{z-z_k(t)}+O(1),
\quad z\to z_k(t),
\end{equation*}
where
\begin{equation}\label{gammak}
\gamma_k(t):=\frac{\dot{A}_k(t)}{D_k(t)}\,.\end{equation}

At the point $z_0(t)$, the function $\dot{f}(z,t)$ has a pole of
order at most $2$, and $f'(z,t)$ has a pole of order $2$. Thus,
$h(z,t)$ has a removable singularity at the point. In more
details, denoting by $d_{-1}(t)$ the residue of $f(z,t)$ at the
point $z_0(t)$, we have
\begin{equation*}\label{z0}
f(z,t)=\frac{d_{-1}(t)}{z-z_0(t)}\,+d_0(t)+O(1),
\end{equation*}
\begin{equation}\label{z0t}
\dot{f}(z,t)=\dot{z}_0(t)\frac{d_{-1}(t)}{(z-z_0(t))^2}\,+\,\frac{\dot{d}_{-1}(t)}{z-z_0(t)}\,+O(1),
\end{equation}
\begin{equation*}\label{z0z}
f'(z,t)=-\frac{d_{-1}(t)}{(z-z_0(t))^2}\,+O(1).
\end{equation*}
From this we see that in a neighborhood of $z_0(t)$ the function
$h(z,t)$ has the expansion
\begin{equation}\label{hzo}
h(z,t)=-\dot{z}_0(t)+o(1), \quad z\to z_0(t).
\end{equation}
In a similar way, we show that $h(z,t)$ has a removable
singularity at the point~$\overline{z}_0(t)$.

The function $$
F(z,t):=h(z,t)-\sum_{j=1}^4\gamma_j(t)\zeta(z-z_j(t))$$ has only
removable singularities at the points $z_k(t)$, $1\le k\le 4$,
${z}_0(t)$, and $\overline{z}_0(t)$, and at points equivalent to
them (by mod of the lattice).  At other points of the plane it is
holomorphic. Consequently,  it can be extended holomorphically to
the whole plane $\mathbb{C}$.

From (\ref{periodh}) we obtain
\begin{equation}\label{periodg}
F(z+\omega_k(t),t)-F(z,t)=-\dot{\omega}_k(t)-\eta_k(t)\sum_{j=1}^4\gamma_j(t),
\quad k=1,2.
\end{equation}
By (\ref{periodg}), the function  $F$ grows not faster than a
linear function, therefore, $F(z,t)=\alpha(t)z+\beta(t)$. So we
have
\begin{equation}\label{h0}
h(z,t)=\sum_{j=1}^4\gamma_j(t)\zeta(z-z_j(t))+\alpha(t)z+\beta(t).
\end{equation}
From (\ref{hzo}) we find
\begin{equation}\label{beta0}
\beta(t)=-\sum_{j=1}^4\gamma_j(t)\zeta(z_0(t)-z_j(t))-\alpha(t)z_0(t)-\dot{z}_0(t).
\end{equation}
From (\ref{periodg}) it follows that
\begin{equation}\label{periods}
\alpha(t)\omega_k(t)=-\dot{\omega}_k(t)-\eta_k(t)\sum_{j=1}^4\gamma_j(t),
\quad k=1,2.
\end{equation} If we put $k=1$, then, taking into account that $\omega_1(t)\equiv
1$, we obtain
\begin{equation}\label{alpha0}
\alpha(t)=-\eta_1(t)\sum_{j=1}^4\gamma_j(t).
\end{equation}
At last, from (\ref{h0}),  (\ref{beta0}), and  (\ref{alpha0}) we
deduce that
\begin{equation}\label{h1}
h(z,t)=\sum_{j=1}^4\gamma_j(t)[\zeta(z-z_j(t))-\zeta(z_0(t)-z_j(t))-\eta_1(t)(z-z_0(t))]-\dot{z}_0(t).
\end{equation}
If we put $k=2$, from  (\ref{periods}) we have
$$
\dot{\omega}_2(t)=-\alpha(t)\omega_2(t)-\eta_2(t)\sum_{j=1}^4\gamma_j(t)=(\omega_2(t)\eta_1(t)-\eta_2(t))\sum_{j=1}^4\gamma_j(t),
$$
and, with the help of the equality (\ref{rel}), we obtain
\begin{equation}\label{period2}
\dot{\omega}_2(t)=2\pi i\sum_{j=1}^4\gamma_j(t).
\end{equation}

Therefore, we proved the following result.

\begin{theorem}\label{loew-ell}
The family $f(z,t)$ satisfies the PDE
$$ \frac{\dot{f}(z,t)}{f'(z,t)}=h(z,t)$$ where $h(z,t)$ is defined by (\ref{h1}); here  $\gamma_k(t)$ and $D_k(t)$ are specified by (\ref{gammak})  and (\ref{dk}). The period $\omega_1(t)$
is equal $1$  and the period $\omega_2(t)$
satisfies~(\ref{period2}).
\end{theorem}

 Now we will write a system of differential equations to find $z_l(t)$, $1\le
l\le 4$. For this, we will write $\dot{f}'(a_l(t),t)$ in two
different ways. On the one hand, from (\ref{fprime}) it follows
that
\begin{equation}\label{aldotprime1}
\dot{f}'(z_l(t),t)=-\dot{z}_l(t)D_l(t).
\end{equation}
On the other hand, by Theorem~\ref{loew-ell}, we have
$\dot{f}(z,t)=h(z,t)f'(z,t)$, therefore,
\begin{multline}\label{dotf}
\dot{f}(z,t)=c(t)\left[\sum_{j=1}^4\gamma_j(t)[\zeta(z-z_j(t))-\zeta(z_0(t)-z_j(t))-
\eta_1(t)(z-z_0(t))]-\dot{z}_0(t)\right]\, \\
\times e^{\gamma(t)
z}\frac{\prod_{k=1}^4\sigma(z-z_k(t))}{\sigma^2(z-z_0(t))\sigma^2(z-\overline{z_0(t)})}
\end{multline}
 and
\begin{multline}\label{dotfpr}
\dot{f}'(z,t)=c(t)\Biggl\{\Biggl[\sum_{j=1}^4\gamma_j(t)\,[\zeta(z-z_j(t))-\zeta(z_0(t)-z_j(t))-\eta_1(t)(z-z_0(t))]-\dot{z}_0(t)\Biggr]\\
\times
\Bigl(\gamma(t)+\sum_{s=1}^4\zeta(z-z_s(t))-2\zeta(z-z_0(t))-2\zeta(z-\overline{z_0(t)})\Bigr)\\
-
\sum_{j=1}^4\gamma_j(t)[\mathfrak{P}(z-z_j(t))-\eta_1(t)]\Biggr\}\,
\,e^{\gamma(t)
z}\frac{\prod_{k=1}^4\sigma(z-z_k(t))}{\sigma^2(z-z_0(t))\sigma^2(z-\overline{z_0(t)})}.
\end{multline}
From (\ref{dotfpr}) we obtain, as $z\to z_l(t)$,
\begin{multline}\label{aldotprime2}
\dot{f}'(z_l(t),t)=c(t)\Biggl[
-\dot{z}_0(t)+\sum_{j=1, j\neq
l}^4\gamma_j(t)\,\bigl[\zeta(z_l(t)-z_j(t))-\zeta(z_0(t)-z_j(t))
\\
-\eta_1(t)(z_l(t)-z_0(t))\bigr] +
\gamma_l(t)\,\Bigl(\sum_{s=1,s\neq l}^4\zeta(z_l(t)-z_s(t))
+\gamma(t)-\eta_1 (z_l(t)-z_0(t))
\\ -\zeta(z_l(t)-z_0(t))-2\zeta(z_l(t)-\overline{z_0(t)})\Bigr)\Biggr]
\,\frac{e^{\gamma(t) z_l(t)}\prod_{k\neq
l,k=1}^4\sigma(z_l(t)-z_k(t))}{\sigma^2(z_l(t)-z_0(t))\sigma^2(z_l(t)-\overline{z_0(t)})}.
\end{multline}

Comparing (\ref{aldotprime1}) and  (\ref{aldotprime2}), taking
into account (\ref{dk}), we see that
\begin{multline}\label{al}
\dot{z}_l=\dot{z}_0- \sum_{j=1, j\neq
l}^4\gamma_j\,\bigl[\zeta(z_l-z_j)-\zeta(z_0-z_j)-\eta_1(z_l-z_0)\bigr]
\\- \gamma_l\,\Bigl(\sum_{s=1,s\neq
l}^4\zeta(z_l-z_s) +\gamma-\eta_1
(z_l-z_0)-\zeta(z_l-z_0)-2\zeta(z_l-\overline{z}_0)\Bigr),\
1\le l \le n.
\end{multline}

Now we will find a differential equation to determine $c(t)$.
Comparing (\ref{family}), (\ref{dotf}), and (\ref{z0t}), we have
\begin{equation}\label{d}
 d_{-1}(t)=-c(t)e^{\gamma(t)
z_0(t)}\frac{\prod_{k=1}^4\sigma(z_0(t)-z_k(t))}{\sigma^2(z_0(t)-\overline{z_0(t)})}\,,
\end{equation}
$$
\dot{d}_{-1}(t)=-c(t)\left\{\sum_{j=1}^4\gamma_j(t)\mathfrak{P}(z_0(t)-z_j(t))+\eta_1(t)+\dot{z}_0(t)\right.
$$
$$
\times\left. \left[\gamma(t)+\sum_{k=1}^4\zeta(z_0(t)-z_k(t))-2\zeta(z_0(t)-\overline{z}_0(t))\right]\right\}\\
\,e^{\gamma(t)
z_0(t)}\frac{\prod_{k=1}^4\sigma(z_0(t)-z_k(t))}{\sigma^2(z_0(t)-\overline{z_0(t)})}\,.
$$
Since
\begin{equation}\label{ident}
\gamma +\sum_{k=1}^4\zeta(z_0-z_k)-2\zeta(z_0-\overline{z}_0)=0,
\end{equation}
we have
\begin{equation*}\label{ddt}
\dot{d}_{-1}(t)=-c(t)\left[\sum_{j=1}^4\gamma_j(t)\mathfrak{P}(z_0(t)-z_j(t))+\eta_1(t)\right]
e^{\gamma(t)
z_0(t)}\frac{\prod_{k=1}^4\sigma(z_0(t)-z_k(t))}{\sigma^2(z_0(t)-\overline{z_0(t)})}\,.
\end{equation*}

Therefore,
\begin{equation}\label{a}
\dot{a}(t)=\sum_{j=1}^4\gamma_j(t)\mathfrak{P}(z_0(t)-z_j(t))+\eta_1(t)
\end{equation}
where $a=\log d_{-1}$.

Differentiating (\ref{ident}), we obtain
$$
i\frac{4\beta}{\pi}\,\frac{\partial \zeta(1/2)}{\partial
\omega_2}\,\dot{m}-\sum_{k=1}^4\mathfrak{P}(z_0-z_k)(\dot{z}_0-\dot{z}_k)+
i2\sum_{k=1}^4\frac{\partial \zeta(z_0-z_k)}{\partial
\omega_2}\,\dot{m}$$
$$+2\mathfrak{P}(z_0-\overline{z}_0)(\dot{z}_0-\dot{\overline{z}}_0)-i4\frac{\partial
\zeta(z_0-\overline{z}_0)}{\partial \omega_2}\,\dot{m}=0,
$$
$$
\left(4\mathfrak{P}(z_0-\overline{z}_0)-\sum_{k=1}^4\mathfrak{P}(z_0-z_k)\right)\dot{z}_0=$$
$$-\sum_{k=1}^4\mathfrak{P}(z_0-z_k)\dot{z}_k+i\left[4\,\frac{\partial
\zeta(z_0-\overline{z}_0)}{\partial
\omega_2}\,-\frac{4\beta}{\pi}\,\frac{\partial
\zeta(1/2)}{\partial \omega_2}\,- 2\sum_{k=1}^4\frac{\partial
\zeta(z_0-z_k)}{\partial \omega_2}\right]\dot{m},
$$
and, therefore,
\begin{multline}\label{yy}
\dot{y}_0=-\sum_{k=1}^4\,\Im\frac{\mathfrak{P}(z_0-z_k)}{4\mathfrak{P}(z_0-\overline{z}_0)-\sum_{j=1}^4\mathfrak{P}(z_0-z_j)}\,\dot{x}_k
+\Re\Biggl[\frac{4\,\displaystyle\frac{\partial
\zeta(z_0-\overline{z}_0)}{\partial
\omega_2}\,}{{4\mathfrak{P}(z_0-\overline{z}_0)-\sum_{k=1}^4\mathfrak{P}(z_0-z_k)}}\\
+\frac{-\displaystyle\frac{4\beta}{\pi}\,\displaystyle\frac{\partial
\zeta(1/2)}{\partial \omega_2}-
2\sum_{k=1}^4\displaystyle\frac{\partial \zeta(z_0-z_k)}{\partial
\omega_2}-\mathfrak{P}(z_0-z_3)+\mathfrak{P}(z_0-z_4)}{{4\mathfrak{P}(z_0-\overline{z}_0)-\sum_{k=1}^4\mathfrak{P}(z_0-z_k)}}\,\Biggr]\dot{m}.
\end{multline}

\begin{theorem}\label{system}
The accessory parameters satisfy the system of ODEs: (\ref{al}),
(\ref{a}), and  (\ref{yy}) where $a=\log d_{-1}$ and $d_{-1}$ is
 defined by (\ref{d}).
\end{theorem}

\begin{corollary}\label{mod}
The conformal module of the domains satisfies the equation
\begin{equation*}
\dot{m}(t)=\pi \sum_{j=1}^4\gamma_j(t).
\end{equation*}
where $\gamma_k(t):={\dot{A}_k(t)}/{D_k(t)}$, $D_k(t)=f''(z_k)$.
\end{corollary}

\section{Symmetric case. Numeric results}\label{numeric}

Now we will describe an approximate method  of finding the
accessory parameters in (\ref{f}). It is based on
Theorem~\ref{system}. If we consider a smooth one-parametric
family $f(z,t)$, $0\le t\le 1$, of conformal mappings of the form
(\ref{family}), then,  knowing the values of the parameters for
$t=0$, we can solve the Cauchy problem with this initial data and
obtain the values of the accessory parameters for all $t$. We note
that it is natural to use the uniform motion of the points
$A_k=A_k(t)$, therefore, in our calculations we will take
$\dot{A}_k=\text{const}$. Moreover, if we choose the appropriate
initial data, then  we change only two of $A_k$, say, $A_1$ and
$A_2$; thus, we can put $\dot{A}_3=\dot{A}_4\equiv0$.

Therefore, to solve the Cauchy problem for the obtained system, we
need to know the initial data, i.e. the values of the accessory
parameters for some $t$. For this, it is convenient to use the
data for the symmetric case when the segment $A_1A_2$ and $A_3A_4$
are symmetric with respect to the real axis and the straight
lines, containing these segments, pass through the origin. (This
can be achieved by a rotation and a shift.)

Now we describe the conformal mapping for the symmetric case.
Because of the Riemann-Schwarz symmetry principle, we can consider
the conformal mapping of a strip onto the upper half of the
symmetric domain $G=G(A_1,A_2,A_3,A_4)$ and then extend it up to
the conformal mapping of the  strip, with twice the original
width, onto the whole domain~$G$.

a) If $0<\beta<\pi$, then the conformal mapping has the form (see
\cite[Part~B, Section~8.2, Example~1, p.~354]{kop_sht}):
\begin{equation*}\label{fsym_theta}
f(z)=\widetilde{c}\,\frac{\vartheta_1(z-\alpha)}{\vartheta_1(z+\alpha)}\,,
\quad \alpha=\frac{\beta}{4\pi}\,,
\end{equation*}
where $\vartheta_1(z)$ is the Jacobi theta-function defined by
(\ref{theta1}) and  $\widetilde{c}>0$ is a constant. From
(\ref{connect}), taking into account that $\omega_1=1$, we easily
deduce that
\begin{equation}\label{fsym}
f(z)=ce^{2\alpha\eta_1z}\,\frac{\sigma(z-\alpha)}{\sigma(z+\alpha)}\,,\quad
c>0.
\end{equation}
We should note that, in contrast to (\ref{f}), here $\sigma(z)$,
defined by (\ref{sw}), matches to the periods $1$ and $im$, not to
$1$ and $i2m$.

\begin{figure}[ht] \centering
\includegraphics[width=4.5 in,%
]{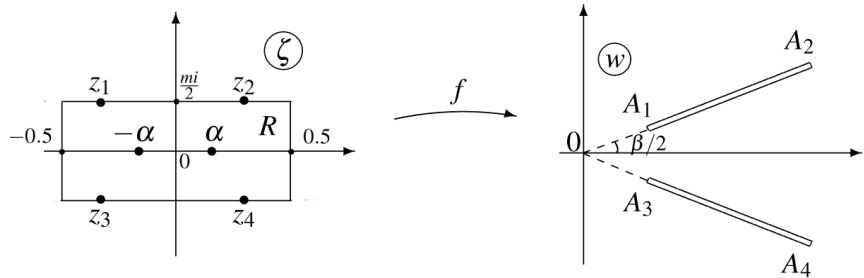}\ \caption{Conformal mapping of the rectangle
$\widetilde{R}$ with identified vertical sides  onto symmetric
domain $G(A_1,A_2,A_3,A_4)$.} \label{symcase}
\end{figure}

The function $f$, defined by~(\ref{fsym}), maps the rectangle
$R:=\{-1/2<\Re z<1/2,\ 0<\Im z<m/2\}$ with identified vertical
sides onto the upper half of $G(A_1,A_2,A_3,A_4)$ and keeps the
real axis; it can be extended, by symmetry, to the rectangle
$\widetilde{R}:=\{-1/2<\Re z<1/2,\ -m/2<\Im z<m/2\}$, and the
extended function maps $\widetilde{R}$ onto the whole domain
$G(A_1,A_2,A_3,A_4)$. The function $f$ has four critical points
$\pm z_k$, $1\le k \le 4$, and $z_1=x_1+im/2$, $z_2=x_2+im/2$,
$z_3=x_1-im/2$, $z_4=x_2-im/2$. Besides, $f$ has a pole at the
point $z=-\alpha$ and a zero at $z=\alpha$ (Fig.~\ref{symcase}).

The critical points $z_k$, $1\le k\le 4$,  can be found from the
equality $f'(z)=0$, i.e.
$$
\zeta(\alpha-z)+\zeta(\alpha+z)=2\eta_1\alpha.
$$
Because of the equality (\cite{akhiezer}, ch.III, \S~15),
\begin{equation*}
\zeta(u+v)+\zeta(u-v)-2\zeta(u)=\frac{\mathfrak{P}'(u)}{\mathfrak{P}(u)-\mathfrak{P}(v)}\,,
\end{equation*}
we have \begin{equation}\label{crit}
\mathfrak{P}(z)=\mathfrak{P}(\alpha)-\frac{\mathfrak{P}'(\alpha)}{2(\alpha\eta_1-\zeta(\alpha))}\,,
\end{equation}
therefore, $z_k$ can be found via the inverse function
$\mathfrak{P}^{-1}$. Because of the evenness of the
$\mathfrak{P}$-function, we see that $z_3=-z_2$ and $z_4=-z_1$.

Without loss of generality we can assume that the nearest points
of the slits are located at the distance $1$ from the origin. Then
the farthest points are at the distance
$l:=if(z_3)/f(z_2)=if(-z_2)/f(z_2)$. Therefore, making use of
(\ref{fsym}) and oddness of the $\sigma$-function, we have
\begin{equation}\label{leng}
l=ie^{-4\alpha\eta_1z_2}\,\frac{\sigma^2(z_2+\alpha)}{\sigma^2(z_2-\alpha)}\,.
\end{equation}

Let $z_2$ be a root of (\ref{crit}); we note that it  depends on
$m$. Then we solve (\ref{leng}) with respect to $m$, to obtain the
initial value of the module. After that, we easily find the
initial values of $z_k$, $1\le k\le 4$. To use them in the
non-symmetric case, we need to shift the obtained values of $z_k$
by the vector $\alpha-i m/2$.

To find $c$ we use the equalities
$$
f(z_3)=f(-z_2)=ce^{-2\alpha\eta_1z_2}\,\frac{\sigma(z_2+\alpha)}{\sigma(z_2-\alpha)}\,,\quad
f(z_2)=ce^{2\alpha\eta_1z_2}\,\frac{\sigma(z_2-\alpha)}{\sigma(z_2+\alpha)}\,.
$$
Multiplying them, we have $c^2=f(z_3)f(z_2)=|f(z_3)f(z_2)|$,
therefore,
$$
c=\sqrt{|f(z_3)f(z_2)|}=\sqrt{|f(z_3)f(z_4)|}.
$$

The residue of $f(z)$, defined by (\ref{fsym}), is equal to
$$
d_{-1}^{\,\,0}=-ce^{-2\alpha^2\eta_1}{\sigma(2\alpha)}, \quad
$$
therefore, the initial value of $a$ is
$$
a^0=\log d_{-1}^{\,\,0}=(1/2)\log|f(z_3)f(z_4)| - 2\alpha^2\eta_1+
\log\sigma (2\alpha)+\pi i.
$$
We note that $|f(z_3)|$ and $|f(z_4)|$ are the distances $l_3$ and
$l_4$ from $A_3$ and $A_4$ to $A_5$;  here $A_5$ is the point of
intersection of the straight lines containing the slits. Finally,
we have
$$
a^0=\log d_{-1}^{\,\,0}=(1/2)\log(l_3l_4)-2\alpha^2\eta_1+
\log\sigma (2\alpha)+\pi i.
$$
\vskip 0.3 cm

b) Consider the case $\beta=0$ when the slits lie on the
(distinct) parallel lines. Without loss of generality we can
assume that the slits are on straight lines parallel to the real
axis. Then the conformal mapping has the form (see \cite[Part~B,
Section~8.1, Example~1, p.~339]{kop_sht} ):
\begin{equation*}\label{fsym1}
f(z)=-\frac{b}{\pi}\,\left(\zeta(z)-\eta_1 z\right).
\end{equation*}
As in the case a), $\zeta(z)$, defined by (\ref{zw}), has the
periods $1$ and $im$, not $1$ and $i2m$. The parameter $b$ means a
half of the vertical distance between the slits. The critical
points $z_k$ of the map can be found from the equation $f'(z)=0$;
it is equivalent to the equality $\mathfrak{P}(z)=-\eta_1$. Using
the evenness of the $\mathfrak{P}$-function, we see that
$x_1=-x_2$, $z_3=-z_2$ and $z_4=-z_1$. Finding $z_2$ and using the
oddness of $f(z)$, we obtain
$$-\frac{b}{\pi}\,\left(\zeta(z_2)-\eta_1 z_2\right)=l/2$$
where $l$ is the length of each slit. From the last equality we
find the initial value of $m$ and $z_k$. As in the case a), to use
the obtained values, we need to shift them; taking into account
that here $\alpha=0$, we see that the shift parameter is the vector $-im/2$. The
residue of $f(z)$ at $z=0$ equals $-b/\pi$, thus, the initial
value of $a$ is $\log(b/\pi)+\pi i$ or
$$
a^0=\log(\Im(z_1-z_3)/(2\pi))+\pi i.
$$

Now we give the Mathematica code, with commentaries, to calculate
the values of parameters, the module and the capacity of ring
domains with the exterior of two rectilinear slits. For convenience, we divide it into 5 steps.

If $0<\beta<\pi$, we first find the point $A_5$ which is the
intersection of the straight lines containing the slits. We will
assume that $A_3A_4$ does not contain $A_5$; in the opposite case
we renumber the points and use the reflection with respect to the
real axis which, in fact, does not change the desired
parameters. We also assume that $A_4$ is farther from $A_5$ than
$A_3$. If $A_1A_2$ also does not contain $A_5$, then we number the
points so that $A_2$ is farther from $A_5$ than~$A_1$. If $A_1A_2$
contains $A_5$, then we consider that $\arg(A2-A5)/(A4-A5)=\beta$.
In the case, either $A_1=A_5$ or $\arg(A1-A5)/(A4-A5)=\beta\pm
\pi$. Dependence on $t$ describes the uniform movement of points
$A_1(t)$ and $A_2(t)$ along the corresponding segments; $A_3(t)$
and $A_4(t)$  are herewith constant. Therefore, $\dot{A}_k(t)$ are
constant, moreover, $\dot{A}_3(t)=\dot{A}_4(t)=0$.
\medskip

\noindent \textbf{Step~1.} Input of location of the points $A_k$,
$1\le k\le 4$. (Here we take $A_1=-2i$, $A_2=3i$, $A_3=1$,
$A_4=3$.) Finding $A_5$, $\beta$,  $\dot{A}_1$,  and $\dot{A}_2$.
\small
\begin{verbatim}
A1=-2.*I; A2=3.*I; A3=1.; A4=3.;
A5=A2+(A1-A2)*Im[(A4-A2)Conjugate[(A3-A4)]]/
Im[(A1-A2)Conjugate[(A3-A4)]];
l1=Sign[Re[(A1-A5)Exp[-I*beta/2]]]*Abs[A1-A5]; l2=Abs[A2-A5];
l3=Abs[A3-A5]; l4=Abs[A4-A5]; beta=Arg[(A2-A1)/(A4-A3)];
alpha=beta/(4*Pi); Adot1=(l1-l3)Exp[I*beta/2];
Adot2=(l2-l4)Exp[I*beta/2];
\end{verbatim}
\normalsize \noindent \textbf{Step~2.} Defining Weierstrass
elliptic functions  ($\mathfrak{P}(z)$, $\mathfrak{P}'(z)$,
$\zeta(z)$, $\sigma(z)$, $\partial \zeta(z)/\partial\omega_2$)
with periods $\omega_1=1$ and $\omega_2=im$ as functions depending
on complex variable $z$ and $m$. Defining functions $\gamma_k(t)$,
$k=1$, $2$. \small
\begin{verbatim}
wp1[z_,w1_,w2_]:=WeierstrassP[z,WeierstrassInvariants[{w1/2,w2/2}]]; wpp1[z_, w1_,
w2_]:=WeierstrassPPrime[z,WeierstrassInvariants[{w1/2,w2/2}]];
wz1[z_,w1_,w2_]:=WeierstrassZeta[z,WeierstrassInvariants[{w1/2,w2/2}]];
ws1[z_,w1_,w2_]:=WeierstrassSigma[z,WeierstrassInvariants[{w1/2,
w2/2}]]; wi1[x_,y_]:=WeierstrassInvariants[{x,y}]; g2[w1_,w2_]:=
  -4(wp1[w1/2,w1,w2]*wp1[w2/2,w1,w2]+wp1[w1/2,w1,w2]*wp1[(w1+w2)/2,w1,w2]+
  wp1[w2/2,w1,w2]*wp1[(w1+w2)/2,w1,w2]);
wz1primeperiod[z_,w1_,w2_]:=-1/(2*Pi*I)((1/2)wpp1[z,w1,w2]+(wz1[z,w1,w2]-
  z*2*wz1[1/2,w1,w2])wp1[z,w1,w2]+2*wz1[w1/2,w1,w2]*wz1[z,w1,w2]-g2[w1,w2]*z/12);
wp[z_,t_]:=wp1[z,1,2*I*t]; ws[z_,t_]:=ws1[z,1,2*I*t]; wz[z_,t_]:=wz1[z,1,2*I*t];
wzw2[z_,t_]:=wz1primeperiod[z,1,2*I*t]; gamma[t_]:=2*beta/Pi*wz[0.5,m[t]];
gamma1[t_]:=-ws[z0[t]-z1[t],m[t]]*ws[z0[t]-z2[t],m[t]]*ws[z0[t]-z3[t],m[t]]*
  ws[z0[t]-z4[t],m[t]]/(ws[2*z0[t],m[t]])^2*Exp[-(a[t]+gamma[t]*(z1[t]-z0[t]))]*
  (ws[z1[t]-z0[t],m[t]])^2*(ws[z1[t]+z0[t],m[t]])^2/(ws[z1[t]-z2[t],m[t]]*
  ws[z1[t]-z3[t],m[t]]*ws[z1[t]-z4[t],m[t]]);
gamma2[t_]:=-ws[z0[t]-z1[t],m[t]]*ws[z0[t]-z2[t],m[t]]*ws[z0[t]-z3[t],m[t]]*
  ws[z0[t]-z4[t],m[t]]/(ws[2*z0[t],m[t]])^2*Exp[-(a[t]+gamma[t]*(z2[t]-z0[t]))]*
  (ws[z2[t] - z0[t],m[t]])^2*(ws[z2[t]+z0[t],m[t]])^2/(ws[z2[t]-z1[t],m[t]]*
  ws[z2[t]-z3[t],m[t]]*ws[z2[t]-z4[t],m[t]]);
\end{verbatim}
\normalsize \noindent \textbf{Step~3.} Finding initial value of
module, critical points, pole, and constant $a$ . \small
\begin{verbatim}
f1[t_]=wp1[alpha,1,I*t]-wpp1[alpha,1,I*t]/(2(alpha*2*wz1[0.5,1,I*t]-
  wz1[alpha,1,I*t]));
Z1[t_]=InverseWeierstrassP[f1[t],WeierstrassInvariants[{0.5,0.5*I*t}]];
L[t_]=Abs[Exp[-4*alpha*2*wz1[0.5,1,I*t]*Z1[t]]*(ws1[Z1[t]+alpha,1,I*t]/
  ws1[Z1[t]-alpha,1,I*t])^2]-l4/l3; ar=0.1; bl=3.0; Do[m0=(ar+ bl)/2.;
  fc=L[m0]; If[L[bl]*fc>0,bl=m0,ar=m0],{i,70}]; X0=Re[Z1[m0]]; x10=beta/(4*Pi)+X0;
x20=beta/(4*Pi)-X0; x30=beta/(4*Pi)+X0; x40=beta/(4*Pi)-X0;
a20=Pi; y00=m0/2; a10=(1/2)*Log[l3*l4]-(beta/(2*Pi))^2*Re[wz1[0.5,1,I*m0]]+
Log[Abs[ws1[beta/(2*Pi),1,I*m0]]];
\end{verbatim}
\normalsize \noindent \textbf{Step~4.} Solving system of ODEs.
\small
\begin{verbatim}
sol = NDSolve[ {-z1'[t]==Re[Adot1*gamma1[t]*(wz[z1[t]-
z2[t],m[t]]+wz[z1[t]-z3[t],m[t]]+
  wz[z1[t]-z4[t],m[t]]+gamma[t]-2*wz[0.5,m[t]]*(z1[t]-z0[t])-wz[z1[t]-z0[t],
  m[t]]-2*wz[z1[t]+z0[t],m[t]])+Adot2*gamma2[t](wz[z1[t]-z2[t],m[t]]-
  wz[z0[t]-z2[t],m[t]]-2*wz[0.5,m[t]]*(z1[t]-z0[t]))],
-z2'[t]==Re[Adot2*gamma2[t]*(wz[z2[t]-z1[t],m[t]]+wz[z2[t]-z3[t],m[t]]+
  wz[z2[t]-z4[t],m[t]]+gamma[t]-2*wz[0.5,m[t]]*(z2[t]-z0[t])-wz[z2[t]-z0[t],
  m[t]]-2*wz[z2[t]+z0[t],m[t]])+Adot1*gamma1[t] (wz[z2[t]-z1[t],m[t]]-
  wz[z0[t]-z1[t],m[t]]-2*wz[0.5,m[t]]*(z2[t]-z0[t]))],
-z3'[t]==-I*m'[t]+Re[Adot1*gamma1[t](wz[z3[t]-z1[t],m[t]]-wz[z0[t]-z1[t],m[t]]-
  2*wz[0.5,m[t]]*(z3[t]-z0[t]))+Adot2*gamma2[t](wz[z3[t]-z2[t],m[t]]-
  wz[z0[t]-z2[t],m[t]]-2*wz[0.5,m[t]]*(z3[t]-z0[t]))],
-z4'[t]==I*m'[t]+Re[Adot1*gamma1[t](wz[z4[t]-z1[t],m[t]]-wz[z0[t]-z1[t],m[t]]-
  2*wz[0.5,m[t]]*(z4[t]-z0[t]))+Adot2*gamma2[t](wz[z4[t]-z2[t],m[t]]-
  wz[z0[t]-z2[t],m[t]]-2*wz[0.5,m[t]]*(z4[t]-z0[t]))],
m'[t]==Re[Pi*(Adot1*gamma1[t]+Adot2*gamma2[t])],
a'[t]==Adot1*gamma1[t]*wp[z0[t]-z1[t],m[t]]+Adot2*gamma2[t]*wp[z0[t]-z2[t],m[t]]
  +2*wz[0.5,m[t]]*(Adot1*gamma1[t]+Adot2*gamma2[t]),
z0'[t]==I*Im[(4*wp[2*z0[t],m[t]]-wp[z0[t]-z1[t],m[t]]-wp[z0[t]-z2[t],m[t]]-
  wp[z0[t]-z3[t],m[t]]-wp[z0[t] - z4[t],m[t]])^(-1)(-wp[z0[t]-z1[t],m[t]]z1'[t]-
  wp[z0[t]-z2[t],m[t]]z2'[t]-wp[z0[t]-z3[t],m[t]]z3'[t]-wp[z0[t]-z4[t],m[t]]
  z4'[t])]+I*Re[(4*wp[2*z0[t],m[t]]-wp[z0[t]-z1[t],m[t]]-wp[z0[t]-z2[t],m[t]]-
  wp[z0[t]-z3[t],m[t]]-wp[z0[t]-z4[t],m[t]])^(-1)(4*wzw2[2*z0[t],m[t]]
  -(4*beta/Pi)wzw2[0.5,m[t]]-2(wzw2[z0[t]-z1[t],m[t]]+wzw2[z0[t]-z2[t],m[t]]+
  wzw2[z0[t]-z3[t],m[t]]+wzw2[z0[t]-z4[t],m[t]]))]*Pi*Re[(Adot1*gamma1[t]+
  Adot2*gamma2[t])],
z1[0]==x10,z2[0]==x20,z3[0]==x30+I*m0,z4[0]==x40-I*m0,a[0]==a10+I*a20,
m[0]==m0,z0[0]==-I*y00},{z1,z2,z3,z4,m,a,z0},{t,0,1.}];
\end{verbatim}
\normalsize \noindent \textbf{Step~5.} Output of desired values of
capacity, module, critical points, pole, and constant~$a$. \small
\begin{verbatim}
s=1.; {1/m[s], m[s], z1[s], z2[s], z3[s], z4[s], a[s], z0[s]} /.sol
\end{verbatim}
\normalsize

In the case of slits, parallel to the real axis, we have the same
system of ODEs. We recall that we can assume that $\Re A_1<\Re
A_2$, $\Re A_3<\Re A_4$, and $\Im A_1=\Im A_2>\Im A_3=\Im A_4$.
Then we find the values of $\dot{A}_1$  and $\dot{A}_2$ by the
formulas $\dot{A}_1=\Re(A_1-A_3)$. $\dot{A}_2=\Re(A_2-A_4)$. We
also have other formulas to find the initial values. Thus, Steps 1
and 3 must be changed to the following ones. \medskip

\noindent \textbf{Step~1'.} Input of location of the points $A_k$,
$1\le k\le 4$. (Here we take $A_1=i$, $A_2=2+i$, $A_3=-2-i$,
$A_4=-1-i$.) Finding $A_5$, $\beta$,  $\dot{A}_1$,  and
$\dot{A}_2$. \small
\begin{verbatim}
A1=1.*I; A2=2.+1.*I; A3=-2.-1.*I; A4=-1.-1.*I;  Adot1=Re[A1-A3];
Adot2=Re[A2-A4]; beta=0.;
\end{verbatim}\medskip
\normalsize \noindent \textbf{Step~3'.} Finding initial value of
module, critical points, pole, and constant $a$. \small
\begin{verbatim}
g[m_]:=-2*WeierstrassZeta[0.5,WeierstrassInvariants[{0.5,0.5*m*I}]];
h[m_]:=Re[InverseWeierstrassP[g[m],WeierstrassInvariants[{0.5,0.5*m*I}]]];
f[m_]:=Re[(2/Pi)(WeierstrassZeta[h[m]+0.5*m*I,WeierstrassInvariants[{0.5,0.5*m*I}
   ]]-2(h[m]+0.5*m*I)WeierstrassZeta[0.5,WeierstrassInvariants[{0.5,0.5*m*I}]])]-
   2*Abs[A3-A4]/Abs[Im[(A3-A1)]]; ar=0.1; bl=3.; Do[m0=(ar+bl)/2.;
   fc=f[m0]; If[f[bl]*fc>0,bl=m0,ar=m0],{i,70}]; X0=Re[h[m0]]; x10=X0;
x20=-X0; x30=X0; x40=-X0; y00=m0/2; a10=Log[Im[A1-A3]/(2*Pi)]; a20=Pi;
\end{verbatim}
\normalsize

\medskip

\textbf{Example~1.} Consider  the case when the endpoints of one
of the segments are the points $a-0.5$, $a+0.5$ on the real axis and the
endpoints of  the other one  are the points $-i$, $-2i$ of the
imaginary axis. Then $\beta=\pi/2$ and $\alpha=0.125$.

As an initial situation, we take the symmetric case when $a=1.5$. With the
help of (\ref{leng}) and (\ref{crit}) we find the initial module
and the real parts of the critical points:
$$m^0=0.67578477\ldots, \quad \widetilde{x}_2^{\,\,0} =-\widetilde{x}_1^{\,\,0} = 0.22367571\ldots$$

Since in the non-symmetric case we have $\Im z_0=0$ in (\ref{f}),
we use a shift $z\mapsto z+\alpha$ in  the $z$-plane and take, as
an initial data, the values
$x^0_k=\widetilde{x}{\,}^{\,0}_k+\alpha$. Therefore,
$$x_1^0 =x_3^0 = 0.34867571\ldots \quad x_2^0=x_4^0=-0.09867571\ldots$$
Because of symmetry, we have $y^0_0=m^0/2$. Consequently,
$$
z^0_{1}=0.34867571\ldots,\phantom{+im^0} \quad
z^0_{2}=-0.09867571\ldots,\phantom{+im^0}
$$
$$
z^0_{3}=0.34867571\ldots+im^0,\quad
z^0_{4}=-0.09867571\ldots-im^0.
$$

At last $a^0=\log d_{-1}^{\,0}$ where $d_{-1}^{\,0}$ is the
residue of (\ref{f}) at the point $z_0$.

\begin{figure}[ht]
\centering
\includegraphics[width=3.5 in,%
]{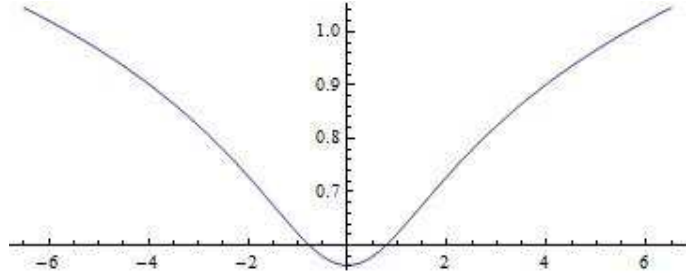}\  \caption{The graph of the dependence of the module $m$
on the parameter $a$ (Example~1).} \label{graph0}
\end{figure}

Finding the residue in the symmetric case, we obtain
$$
a^0=[\ln
\sqrt{2}-4\alpha^2\zeta(0.5;1,im^0)+\ln{\sigma(2\alpha;1,im^0)}]+\pi
i
$$
$$
=-1.11526111\ldots + i\,3.14159265\ldots$$

Here the functions $\zeta(z)=\zeta(z;1,im^0)$ and
$\sigma(z)=\sigma(z;1,im^0)$ correspond to the periods $1$ and
$im^0$. Solving the system of differential equations, we find the
dependence of the parameters in (\ref{f}) on the parameter $a$
(see Fig.~\ref{graph0}).

The values of moduli for some $a$ are given on the
Table~\ref{tab1}.
\medskip

\begin{table}[ht]
\caption{The values of moduli and capacities for some $a$
(Example~1).} \centering
\begin{tabular}
{|c|c|c|c|c|c|c|c|c|
}
  \hline
$a$ &  0 & 1 & 2 & 3 & 4 & 5 & 6 & 7
\\
  \hline
 $m$ &  0.56247 & 0.62207 & 0.72955 & 0.82469 & 0.90239 & 0.96656 & 1.02073 & 1.06743 
 \\
  \hline
 cap &  1.77787& 1.60753& 1.37070 & 1.21258 & 1.10817 & 1.03459 & 0.97968 &0.93682
 \\
  \hline
\end{tabular}\label{tab1}
\end{table}

\textbf{Example~2.} We also computed the moduli $\text{mod}\,G$
and the corresponding capacities $\text{cap}\,G$ for some domains
$G(A_1,A_2,A_3,A_4)$ when $A_k$ are from the integer lattice in
the complex plane. Comparison our results with those obtained by
other methods show very good coincidence, up to $10^{-6}$. In
Table~\ref{tab2} we give some values of capacities obtained by our
method and by a MATLAB algorithm written by Prof. M.~Nasser
\cite{nv}; the values are given with $8$ digits after the decimal
point.\medskip

\begin{table}[ht]
\caption{The values of capacities for some domains
$G(A_1,A_2,A_3,A_4)$ (Example~2).}

\centering
\begin{tabular}{|c|c|c|c|c|c|c|}
  \hline
&&&&&\multicolumn{2}{|c|}{$\text{cap}\,G$}\\
\cline{6-7} &\raisebox{1.5ex}[0cm][0cm]{$A_1$}&\raisebox{1.5ex}[0cm][0cm]{$A_2$}&\raisebox{1.5ex}[0cm][0cm]{$A_3$}&\raisebox{1.5ex}[0cm][0cm]{$A_4$}& our results & Nasser's\\
  \hline
  $1$ &   $i$  &  $2+i$  &   $-2-i$  &  $-1-i$   &1.44058466&  1.44058486\\
  \hline
  $2$  &  $i$   & $2+i$   &  $-2-2i$  &  $-1-2i$ &1.30971558&  1.30971579\\
\hline
  $3$  &  $i$  &  $2+i$   &  $3-2i$  &  $4-3i$   &1.35832035 & 1.35832051\\
\hline
  $4$   & $i$  &  $2+2i$  &  $-2-i$  &  $-1-i$   &1.42710109 & 1.42710150\\
\hline
  $5$ &    $i$  &  $2+2i$  &  $-2-2i$ &  $-1-2i$ &1.29776864 & 1.29776889\\
\hline
  $6$  &  $i$   & $2+2i$   & $3-2i$  &  $4-3i$   &1.32814214 & 1.32814249\\
\hline
  $7$   & $i$  &  $3+2i$   & $-2-i$  &  $-1-i$   &1.49363842 & 1.49363897\\
\hline
  $8$  &    $i$ &   $3+2i$  &  $-2-2i$ &  $-1-2i$&1.36333122&  1.36333156\\
\hline
  $9$  &  $i$  &  $3+2i$  &  $3-2i$  &  $4-3i$   &1.45844055 & 1.45844094\\
\hline
  $10$  &  $i$  &  $3i$  &  $3$  &  $4$          &1.29126199 & 1.29126229\\
\hline
  $11$  &  $i$  &  $3i$  &  $0$  &  $2$          &2.18251913&  2.18251948\\
\hline
  $12$  &  $i$  &  $3i$  &  $-3$  &  $2$         &2.82846257 & 2.82846345\\
\hline
  $13$  &  $i$  &  $3+i$  &  $-i$  &  $3-i$      &2.69941565 & 2.69941690\\
\hline
  $14$  &  $i$  &  $3+2i$  &  $-i$  &  $3-2i$    &2.23470313&  2.23470399\\
\hline
  $15$  &  $i$  &  $3+3i$  &  $-i$  &  $3-3i$    &2.11547784&  2.11547801\\
  \hline
\end{tabular}\label{tab2}
\end{table}

\section{Monotonicity of conformal module }\label{extr}

Now we will investigate behavior of the conformal module of
$G=G(A_1,A_2,A_3,A_4)$ in the case when the segment $A_3A_4$ is
fixed and the segment $A_1A_2$ slides along a straight line with a
fixed length. This case is equivalent to the situation when the
segment $A_1A_2$ is fixed and $A_3A_4$ has a fixed length and
shifts by vectors with a fixed direction. 

We note that some similar problems for quadrangles were
investigated by Dubinin and Vuorinen~\cite{dub_vuor}.

Without loss of generality we may assume that $A_1A_2$ lies on the
real axis and $A_1$ is the left endpoint of the segment. Moreover,
we can consider the family with $A_1=t$, $A_2=t+l$, where $l$ is
the length of $A_1A_2$. Then  $\dot{A}_1(t)=\dot{A}_2(t)=1$. It is
clear that $\dot{A}_3(t)=\dot{A}_4(t)=0$. From Corollary~\ref{mod}
we obtain that $\dot{m}(t)=\pi(\gamma_1(t)+\gamma_2(t))$ where
$\gamma_k(t)=1/f''(x_k,t)$, $k=1$,~$2$. It is easy to see that
$f''(x_1,t)>0$ and $f''(x_2,t)<0$. Therefore,
$$\dot{m}(t)=\pi\left(|\gamma_1(t)|-|\gamma_2(t)|\right)=\pi\left(\frac{1}{|f''(x_1,t)|}-\frac{1}{|f''(x_2,t)|}\right).$$

If $|f''(x_1,t)|>|f''(x_2,t)|$, then, when moving a segment
$A_1A_2$ to the right, the conformal module of
$G=G(A_1,A_2,A_3,A_4)$ decreases, otherwise, it increases. At
critical points of the module we have $|f''(x_1,t)|=|f''(x_2,t)|$.

Now we compare $|f''(x,t)|$ at the points $x_1$  and $x_2$ using
methods of the symmetrization theory. We will temporarily assume
that $A_1A_2$ is symmetric with respect to the imaginary axis,
i.e. $A_2$ lies on the positive part of the real axis
symmetrically to $A_1$.

In the following lemma we investigate a more general case when the
considered doubly connected domain $G$ is the exterior of the
segment $A_1A_2$ and some continuum $Q$; if $Q=A_3A_4$, we obtain
our case.

\begin{lemma}\label{sym}
Let the continuum $Q$ lie in the right half-plane $\Re w>0$ and
let $\psi: \{q<|\zeta|<1\}\to G$ be  a conformal mapping. If $\zeta_1$
and $\zeta_2$ are the points of the unit circle corresponding to
the endpoints $A_1$ and $A_2$ of the segment, then
$|\psi''(\zeta_1)|>|\psi''(\zeta_2)|$.
\end{lemma}

\begin{figure}[ht]
\centering
\includegraphics[width=3.8 in,%
]{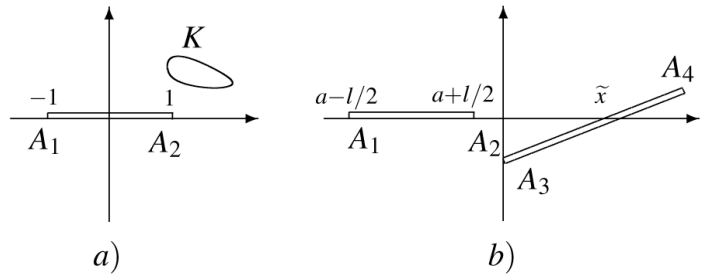}\
\caption{The domain $G$: $a)$  in Lemma~\ref{sym}; \ $b)$ 
in Corollary~\ref{corseg}.} \label{graph3}
\end{figure}

Proof. Without loss of generality we can assume that $A_1A_2$
coincides with the segment $[-1,1]$ (Fig.~\ref{graph3}, $a)$).
Consider the function $\varphi$ inverse to the Joukowskii
function; it maps $G$ onto the exterior of the unit disk with
excluded set $Q_1:=\varphi(Q)$. Using the Riemann-Schwarz symmetry
principle, we conclude that $\varphi(Q)$ lies in the right
half-plane. Now applying the symmetry principle once more, we can
extend $\varphi\circ \psi$ to the annulus $A:=\{q<|\zeta|<1/q\}$.
The function $\varphi\circ \psi$ maps it onto the doubly connected
domain $G_1:=\overline{\mathbb{C}}\setminus (Q_1\cap Q_2)$; here
$Q_2$ is symmetric to $Q_1$ with respect to the unit circle,
therefore, it also  lies in the right half-plane.

Now consider the reduced moduli of $A$ at the points $\zeta_1$ and
$\zeta_2$; it is obvious that they are equal, i.e.
$r(A,\zeta_1)=r(A,\zeta_2)$. On the other side,
$$r(G_1,-1)=r(A,\zeta_1)+\frac{1}{4\pi}\,\log |\psi''(\zeta_1)|,\quad
r(G_1,1)=r(A,\zeta_2)+\frac{1}{4\pi}\,\log |\psi''(\zeta_2)|.$$
Therefore, we only need to show that $r(G_1,-1)>r(G_1,1)$. But this conclusion
 follows from \cite{dub}, thrm.~1.2, because the configuration
$(G_1,-1)$ is obtained from $(G_1,1)$ by polarization. \vskip 0.5
cm

\begin{corollary}\label{corseg}
Let  $A_3A_4$ be a fixed  segment  in the right half-plane,
intersecting the real at the point $\widetilde{x}\,,$
and let one of its
endpoints lie on the imaginary axis. Let $A_1A_2$ be the segment
$[a-l/2,a+l/2]$ on the real axis with a fixed length $l$
(Fig.~\ref{graph3}, $b)$). If $\widetilde{x}\le l/2$, then, when\
\, $a$\ \, increases from $-\infty$ to $\widetilde{x}-l/2$, the
conformal module of $G(A_1,A_2,A_3,A_4)$ decreases from $+\infty$
to $0$. If $\widetilde{x}> l/2$, then the conformal module
decreases from $+\infty$ to some positive value, when \ \, $a$\ \,
increases from $-\infty$ to $0$.
 \end{corollary}

If $A_3A_4$ does not intersect the real axis then the conformal
module decreases for $a$ close to $-\infty$ and increases for $a$
close to $+\infty$. In this connection the problem arises: does the
module always have a unique minimum or are there  situations when it
has more than one (local) minimum?

It is also interesting to investigate the problem for the case
when the slits are parallel to each other. Then, using the result
that the conformal module decreases after symmetrization with
respect to a straight line, we  conclude that the minimum of the
conformal module is attained for the case of slits symmetric with
respect to the orthogonal line.

The same is valid when the slits are perpendicular to each other,
one of the slits is fixed and does not intersect the straight line
containing the second one. Then the minimal module is attainted
for the case when the second slit is symmetric with respect to the
line containing the first slit.



\begin{thebibliography}{bbgghv}

\bibitem{as}
{\sc M.~Abramowitz and I. A. Stegun,} Handbook of mathematical
functions with formulas, graphs, and mathematical tables. National
Bureau of Standards Applied Mathematics Series, 55 For sale by the
Superintendent of Documents, U.S. Government Printing Office,
Washington, D.C. 1964 xiv+1046 pp.

\bibitem{akhiezer}
{\sc N.I.~Akhiezer,}  Elements of the Theory of Elliptic Functions.
Transl. of Mathematical Monographs, vol. 79, American Mathematical
Soc., RI, 1990.

\bibitem{alex} {\sc I.A.~Aleksandrov,} Parametric continuations in the theory of univalent functions.
Nauka, Moscow, 1976 (Russian).

\bibitem{bsv} {\sc D.
Betsakos, K. Samuelsson, and M. Vuorinen, } The computation of capacity of planar condensers. Publ. Inst. Math. (Beograd) (N.S.) 75(89) (2004), 233--252.

\bibitem{bbgghv}{\sc S. Bezrodnykh, A. Bogatyrev, S. Goreinov, O. Grigoriev, H. Hakula,  and M. Vuorinen}:
{On capacity computation for symmetric polygonal condensers}.
{J. Comput. Appl. Math.  361 (2019), 271--282, }  {https://doi.org/10.1016/j.cam.2019.03.030}

\bibitem{gum1} {\sc M.D.~Contreras, S.~Diaz-Madrigal, and P.~Gumenyuk,}
Loewner theory in annulus I: Evolution families and differential
equations. Trans. Amer. Math. Soc. 365 (2013), 2505--2543.

\bibitem{gum2} {\sc M.D.~Contreras, S.~Diaz-Madrigal, and P.~Gumenyuk,}
Loewner theory in annulus II: Loewner chains. Analysis and
Mathematical Physics, 2011, Volume 1, Issue 4, 351--385.

\bibitem{cr} {\sc D.
Crowdy,} Schwarz-Christoffel mappings to unbounded multiply connected polygonal regions. Math. Proc. Cambridge Philos. Soc. 142 (2007), no. 2, 319--339.

\bibitem{delil1} {\sc T.K.~DeLillo, A.R.~Elcrat, and J.A.~Pfaltzgraff,}
Schwarz--Christoffel Mapping of the Annulus. SIAM Rev., 2001,
43(3), 469--477.

\bibitem{delil2}  {\sc T.K.~DeLillo, T.A.~Driscoll, A.R.~Elcrat, and
J.A.~Pfaltzgraff,} Computation of Multiply Connected
Schwarz-Christoffel Maps for Exterior Domains. Computational
Methods and Function Theory. 2006, Volume 6, Issue 2, pp 301--315.

\bibitem{dt} {\sc T.A.~Driscoll and L.N.~Trefethen,} Schwarz-Christoffel mapping. Cambridge Monographs on Applied and Computational Mathematics, 8. Cambridge University Press, Cambridge, 2002. xvi+132 pp.

\bibitem{dub}
{\sc V.N.~Dubinin,} Symmetrization in the geometric theory of functions
of a complex variable. Russian Mathematical Surveys, 1994, 49:1,
1--79.

\bibitem{dubook}  \textsc{V.N.~Dubinin,} {Condenser capacities and symmetrization in geometric function theory}. {Translated from the Russian by Nikolai G. Kruzhilin},
 {Springer, Basel}, {2014}, {xii+344}.

\bibitem{dub_vuor} {\sc V.N.~Dubinin and  M.~Vuorinen,}  On conformal moduli of polygonal quadrilaterals. Israel J. Math. 171,
111--125 (2009).

\bibitem{garnett}
{\sc J.B.~Garnett,  D.E.~Marshall,} Harmonic measure. Reprint of
the 2005 original. New Mathematical Monographs, 2. Cambridge
University Press, Cambridge, 2008. xvi+571 pp.

\bibitem{gol1} {\sc G.M.~Goluzin,} On the parametric representation of functions
univalent in a ring. Mat. Sb. (N.S.), 29(71):2 (1951), 469--476
(Russian).

\bibitem{gol} {\sc G.M.~Goluzin,} Geometric Theory of Functions of a Complex
Variable. Translations of Mathematical Monographs, AMS,  1969.

\bibitem{hrv1} { \textsc{H.
Hakula, A. Rasila and M. Vuorinen,} On moduli of rings and quadrilaterals: algorithms and experiments. SIAM J. Sci. Comput. 33 (2011), no. 1, 279--302. }

\bibitem{hrv3} {\sc H.~
Hakula,  A.~Rasila, and  M.~Vuorinen,} Conformal modulus and
planar domains with strong singularities and cusps. Electron.
Trans. Numer. Anal. 48 (2018), 462--478.

\bibitem{henrici3} {\sc P.~Henrici,} Applied and Computational Complex Analysis,
Vol.~3: Discrete Fourier Analysis, Cauchy Integrals, Construction
of Conformal Maps, Univalent Functions. Wiley, New York, 1986.

\bibitem{komatu1} {\sc Yu.~Komatu,} Untersuchungen  \"uber konforme
Abbildung zweifach zusammenh\"angender  Bereiche. Proc. Phys.,
Math. Soc. Japan 25(1943), 1--42

\bibitem{komatu} {\sc Yu.~Komatu,} Darstellungen der in einem
Kreisringe analytischen Funktionen nebst den Anwendungen auf
konforme Abbildung \"uber Polygonalringgebiete. Jap. J. Math. 19
(1945), 203-215 (German).

\bibitem{kop_sht} {\sc W.~Koppenfels and F.~Stallmann,} Praxis der konformen Abbildung. Berlin: Springer,
1959 (German).

\bibitem{nv} {\sc M.~Nasser and M.~Vuorinen,} Computation of conformal invariants.
Manuscript, August 2019.

\bibitem{nas_vuz} {\sc S.R.~Nasyrov,}  Uniformization of One-Parametric Families of Complex
Tori. Russian Mathematics, 2017, Vol. 61, No. 8, 36--45.

\bibitem{nas_petr} {\sc S.R.~Nasyrov,}
Families of elliptic functions and uniformization of complex tori
with a unique point over infinity. Probl. Anal. Issues Anal., 2018,
Vol. 7(25):2,  98--111.

\bibitem{nas_dokl} {\sc S.R.~Nasyrov,} Uniformization of Simply-Connected Ramified
Coverings of the Sphere by Rational Functions. Lobachevskii
Journal of Mathematics. 2018,  V.39, No.2 252--258.

\bibitem{nist}{\sc F.W.J.~Olver, D.W.~Lozier, R.F. Boisvert, and Ch.W.Clark} (eds.)
NIST Handbook of Mathematical Functions. U.S. Department of
Commerce, National Institute of Standards and Technology,
Washington, DC; Cambridge University Press, Cambridge, 2010.
xvi+951 pp.

\bibitem{ps} {\sc N. Papamichael and  N.
Stylianopoulos,} Numerical conformal mapping. Domain decomposition
and the mapping of quadrilaterals. World Scientific Publishing Co.
Pte. Ltd., Hackensack, NJ, 2010. xii+229 pp.

\bibitem{dlmf1} {\sc W.P.~Reinhardt, P.L.~Walker} Digital Library of Mathematical functions. Chapter 20. Theta
Functions. https://dlmf.nist.gov/20.

\bibitem{dlmf2} {\sc W.P.~Reinhardt, P.L.~Walker} Digital Library of Mathematical functions. Chapter 23. Weierstrass Elliptic and Modular Functions.
https://dlmf.nist.gov/23.

\end{thebibliography}
\end{document}